\documentclass[dvips,12pt]{article}

\usepackage{etex}
\usepackage{tikz-cd}
\tikzcdset{
  shift left/.default=+0.7ex,% default: +0.56ex
  shift right/.default=+0.7ex}

\newenvironment{tikzar}[1][]{{}\kern-4pt\begin{tikzcd}[ampersand replacement=\&,#1]}%
{\end{tikzcd}\kern-4pt{}}

\usetikzlibrary{cd}
\usetikzlibrary{decorations.pathmorphing,positioning}

\usepackage{amssymb}
\usepackage{amstext}
\usepackage{amsmath,stmaryrd}
\usepackage{txfonts}
\usepackage{bbm}
\usepackage{enumerate}
\usepackage{fullpage}
\usepackage{multicol}
\usepackage{rotating}

\newdimen\proofrulebreadth \proofrulebreadth=.05em
\newdimen\proofdotseparation \proofdotseparation=1.25ex
\newdimen\proofrulebaseline \proofrulebaseline=2ex
\newcount\proofdotnumber \proofdotnumber=3
\let\then\relax
\def\hfi{\hskip0pt plus.0001fil}
\mathchardef\squigto="3A3B
\newif\ifinsideprooftree\insideprooftreefalse
\newif\ifonleftofproofrule\onleftofproofrulefalse
\newif\ifproofdots\proofdotsfalse
\newif\ifdoubleproof\doubleprooffalse
\let\wereinproofbit\relax
\newdimen\shortenproofleft
\newdimen\shortenproofright
\newdimen\proofbelowshift
\newbox\proofabove
\newbox\proofbelow
\newbox\proofrulename
\def\shiftproofbelow{\let\next\relax\afterassignment\setshiftproofbelow\dimen0 }
\def\shiftproofbelowneg{\def\next{\multiply\dimen0 by-1 }%
\afterassignment\setshiftproofbelow\dimen0 }
\def\setshiftproofbelow{\next\proofbelowshift=\dimen0 }
\def\setproofrulebreadth{\proofrulebreadth}

\def\prooftree{% NESTED ZERO (\ifonleftofproofrule)
\ifnum  \lastpenalty=1
\then   \unpenalty
\else   \onleftofproofrulefalse
\fi
\ifonleftofproofrule
\else   \ifinsideprooftree
        \then   \hskip.5em plus1fil
        \fi
\fi
\bgroup% NESTED ONE (\proofbelow, \proofrulename, \proofabove,
\setbox\proofbelow=\hbox{}\setbox\proofrulename=\hbox{}%
\let\justifies\proofover\let\leadsto\proofoverdots\let\Justifies\proofoverdbl
\let\using\proofusing\let\[\prooftree
\ifinsideprooftree\let\]\endprooftree\fi
\proofdotsfalse\doubleprooffalse
\let\thickness\setproofrulebreadth
\let\shiftright\shiftproofbelow \let\shift\shiftproofbelow
\let\shiftleft\shiftproofbelowneg
\let\ifwasinsideprooftree\ifinsideprooftree
\insideprooftreetrue
\setbox\proofabove=\hbox\bgroup$\displaystyle % NESTED TWO
\let\wereinproofbit\prooftree
\shortenproofleft=0pt \shortenproofright=0pt \proofbelowshift=0pt
\onleftofproofruletrue\penalty1
}

\def\eproofbit{% NESTED TWO
\ifx    \wereinproofbit\prooftree
\then   \ifcase \lastpenalty
        \then   \shortenproofright=0pt  % 0: some other object, no indentation
        \or     \unpenalty\hfil         % 1: empty hypotheses, just glue
        \or     \unpenalty\unskip       % 2: just had a tree, remove glue
        \else   \shortenproofright=0pt  % eh?
        \fi
\fi
\global\dimen0=\shortenproofleft
\global\dimen1=\shortenproofright
\global\dimen2=\proofrulebreadth
\global\dimen3=\proofbelowshift
\global\dimen4=\proofdotseparation
\global\count255=\proofdotnumber
$\egroup  % NESTED ONE
\shortenproofleft=\dimen0
\shortenproofright=\dimen1
\proofrulebreadth=\dimen2
\proofbelowshift=\dimen3
\proofdotseparation=\dimen4
\proofdotnumber=\count255
}

\def\proofover{% NESTED TWO
\eproofbit % NESTED ONE
\setbox\proofbelow=\hbox\bgroup % NESTED TWO
\let\wereinproofbit\proofover
$\displaystyle
}%
\def\proofoverdbl{% NESTED TWO
\eproofbit % NESTED ONE
\doubleprooftrue
\setbox\proofbelow=\hbox\bgroup % NESTED TWO
\let\wereinproofbit\proofoverdbl
$\displaystyle
}%
\def\proofoverdots{% NESTED TWO
\eproofbit % NESTED ONE
\proofdotstrue
\setbox\proofbelow=\hbox\bgroup % NESTED TWO
\let\wereinproofbit\proofoverdots
$\displaystyle
}%
\def\proofusing{% NESTED TWO
\eproofbit % NESTED ONE
\setbox\proofrulename=\hbox\bgroup % NESTED TWO
\let\wereinproofbit\proofusing
\kern0.3em$
}

\def\endprooftree{% NESTED TWO
\eproofbit % NESTED ONE
  \dimen5 =0pt% spread of hypotheses
\dimen0=\wd\proofabove \advance\dimen0-\shortenproofleft
\advance\dimen0-\shortenproofright
\dimen1=.5\dimen0 \advance\dimen1-.5\wd\proofbelow
\dimen4=\dimen1
\advance\dimen1\proofbelowshift \advance\dimen4-\proofbelowshift
\ifdim  \dimen1<0pt
\then   \advance\shortenproofleft\dimen1
        \advance\dimen0-\dimen1
        \dimen1=0pt
        \ifdim  \shortenproofleft<0pt
        \then   \setbox\proofabove=\hbox{%
                        \kern-\shortenproofleft\unhbox\proofabove}%
                \shortenproofleft=0pt
        \fi
\fi
\ifdim  \dimen4<0pt
\then   \advance\shortenproofright\dimen4
        \advance\dimen0-\dimen4
        \dimen4=0pt
\fi
\ifdim  \shortenproofright<\wd\proofrulename
\then   \shortenproofright=\wd\proofrulename
\fi
\dimen2=\shortenproofleft \advance\dimen2 by\dimen1
\dimen3=\shortenproofright\advance\dimen3 by\dimen4
\ifproofdots
\then
        \dimen6=\shortenproofleft \advance\dimen6 .5\dimen0
        \setbox1=\vbox to\proofdotseparation{\vss\hbox{$\cdot$}\vss}%
        \setbox0=\hbox{%
                \advance\dimen6-.5\wd1
                \kern\dimen6
                $\vcenter to\proofdotnumber\proofdotseparation
                        {\leaders\box1\vfill}$%
                \unhbox\proofrulename}%
\else   \dimen6=\fontdimen22\the\textfont2 % height of maths axis
        \dimen7=\dimen6
        \advance\dimen6by.5\proofrulebreadth
        \advance\dimen7by-.5\proofrulebreadth
        \setbox0=\hbox{%
                \kern\shortenproofleft
                \ifdoubleproof
                \then   \hbox to\dimen0{%
                        $\mathsurround0pt\mathord=\mkern-6mu%
                        \cleaders\hbox{$\mkern-2mu=\mkern-2mu$}\hfill
                        \mkern-6mu\mathord=$}%
                \else   \vrule height\dimen6 depth-\dimen7 width\dimen0
                \fi
                \unhbox\proofrulename}%
        \ht0=\dimen6 \dp0=-\dimen7
\fi
\let\doll\relax
\ifwasinsideprooftree
\then   \let\VBOX\vbox
\else   \ifmmode\else$\let\doll=$\fi
        \let\VBOX\vcenter
\fi
\VBOX   {\baselineskip\proofrulebaseline \lineskip.2ex
        \expandafter\lineskiplimit\ifproofdots0ex\else-0.6ex\fi
        \hbox   spread\dimen5   {\hfi\unhbox\proofabove\hfi}%
        \hbox{\box0}%
        \hbox   {\kern\dimen2 \box\proofbelow}}\doll%
\global\dimen2=\dimen2
\global\dimen3=\dimen3
\egroup % NESTED ZERO
\ifonleftofproofrule
\then   \shortenproofleft=\dimen2
\fi
\shortenproofright=\dimen3
\onleftofproofrulefalse
\ifinsideprooftree
\then   \hskip.5em plus 1fil \penalty2
\fi
}

\renewcommand{\to}{\xrightarrow{\ \ \ }}%{\longrightarrow}
\newcommand{\tto}[1]{\xrightarrow{#1}}
\newcommand{\oot}[1]{\xleftarrow{#1}}
\newcommand{\mono}{\rightarrowtail}
\newcommand{\epi}{\twoheadrightarrow}

\newcommand{\inclusion}{\hookrightarrow}
\newcommand{\mmono}[1]{\stackrel{#1}\rightarrowtail}
\newcommand{\eepi}[1]{\stackrel{#1}\twoheadrightarrow}

\newcommand{\mapfrom}{\mathrel{\reflectbox{\ensuremath{\mapsto}}}}
\newcommand{\pfn}{\rightharpoonup}

\newcommand{\WP}{\mbox{\large $\wp$}}
\newcommand{\Rel}{\mathsf{Rel}}%{\mathbf{R}}%

\newcommand{\Set}{\mathsf{Set}}%
\newcommand{\set}{\mathsf{set}}%{\mathbf{S}}%{{\mathfrak S}}%
\newcommand{\Cat}{\mathsf{Cat}}

\newcommand{\CatS}{\mathsf{CatS}}
\newcommand{\CatA}{\mathsf{CatA}}
\newcommand{\CaT}{\mathsf{CaT}}

\newcommand{\Esp}{\mathsf{Esp}}%{\mathsf{tp}}
\newcommand{\MLat}{\mathsf{Msl}}
\newcommand{\Ssp}{\mathsf{SobS}}%{\mathsf{tp}}
\newcommand{\Ssl}{\mathsf{SpaL}}

\newcommand{\CaL}{\mathsf{CaL}}
\newcommand{\PSg}{\mathsf{PSg}}

\newcommand{\ccirc}{\circledcirc}

\newcommand{\id}{\iota}%{\mathrm{id}}
\newcommand{\idds}[1]{{#1}_\iota}
\newcommand{\Idds}[1]{{#1}^\iota}
\newcommand{\IIdds}{\Idds{(-)}}

\newcommand{\OOO}{{\cal O}}
\newcommand{\PPP}{{\cal P}}

\newcommand{\CCc}{{\mathbb C}}
\newcommand{\DDd}{{\mathbb D}}

\newcommand{\NNn}{{\mathbb N}}

\newcommand{\PPp}{{\mathbb P}}

\newcommand{\XXx}{{\mathbb X}}
\newcommand{\YYy}{{\mathbb Y}}

\newcommand{\TTwo}{\mathbbm{2}}

\newcommand{\mathbold}[1]{\mbox{\boldmath $#1$}}

\newcommand{\Ccc}{{\mathbold C}}

\mathcode`\<="4268 %left delimiter
\mathcode`\>="5269 %right delimiter
\mathchardef\gt="313E %relation >
\mathchardef\lt="313C %relation <

 \def\pushright#1{{%              set up
    \parfillskip=0pt            % so \par doesnt push \square to left
    \widowpenalty=10000         % so we dont break the page before \square
    \displaywidowpenalty=10000  % ditto
    \finalhyphendemerits=0      % TeXbook exercise 14.32
    \leavevmode                 % \nobreak means lines not pages
    \unskip                     % remove previous space or glue
    \nobreak                    % don't break lines
    \hfil                       % ragged right if we spill over
    \penalty50                  % discouragement to do so
    \hskip.2em                  % ensure some space
    \null                       % anchor following \hfill
    \hfill                      % push \square to right
    {#1}                        % the end-of-proof mark (or whatever)
    \par}}                      % build paragraph

 \def\qed{\pushright{$\square$}\penalty-700 \smallskip}

\newenvironment{prf}[1]{\begin{trivlist} \item[{\bf ~Proof}#1.]}%
{\qed\end{trivlist}}

\newcommand{\beq}{\begin{equation}}
\newcommand{\eeq}{\end{equation}}
\newcommand{\ba}[1]{\begin{array}{#1}}
\newcommand{\ea}{\end{array}}
\newcommand{\bea}{\begin{eqnarray}}
\newcommand{\eea}{\end{eqnarray}}
\newcommand{\bear}{\begin{eqnarray*}}
\newcommand{\eear}{\end{eqnarray*}}
\newcommand{\bpr}{\begin{prf}{}}
\newcommand{\epr}{\end{prf}}
\newcommand{\bprf}[1]{\begin{prf}{#1}}
\newcommand{\eprf}{\end{prf}}

\newtheorem{theorem}{Theorem}[section]
\newtheorem{proposition}[theorem]{Proposition}
\newtheorem{lemma}[theorem]{Lemma}
\newtheorem{corollary}[theorem]{Corollary}

\newtheorem{definition}[theorem]{Definition}

\newcommand{\kar}[1]{{#1}_{\circlearrowleft}}
\newcommand{\Kar}[1]{{#1}^{\circlearrowleft}}
\newcommand{\Karr}{\circlearrowleft}
\newcommand{\Skel}[1]{{#1}^{\boxdot}}
\newcommand{\Skell}{\boxdot}
\newcommand{\Taut}[1]{{#1}^{\boxcircle}}
\newcommand{\Tautt}{\boxcircle}

\newcommand{\tp}{\OOO}%{\mathsf{tp}}
\newcommand{\pt}{\Upsilon}%{\maltese}
\newcommand{\otp}[1]{{#1}_{1}}
\newcommand{\ootp}{\otp{(-)}}
\newcommand{\upt}[1]{{#1}^{\circlearrowleft}}
\newcommand{\uupt}{\upt{(-)}}

\newcommand{\op}{o}

\newcommand{\Sier}{\mathbold{2}}

\newcommand{\halts}{\!\downarrow}

 \newcommand{\wwidetilde}[1]{\widetilde{\widetilde{#1}}}

 \newcommand{\iapprox}{\stackrel\iota\approx}

\begin{document}
\title{Catales: \\
From topologies without points\\
to categories without objects\\[1ex]
%and the sobriety of absolute completions
{\Large(or: How not to be evil categorically)}
}

\author{Dusko Pavlovic\\
(working paper, 2021\ldots)}

\date{}

\maketitle

\begin{abstract}
Locales have been studied as ``topologies without points", mainly by tools of category theory. While traditional topology presents a space as a set of points with specified neighborhoods, localic topology presents a space as a lattice of neighborhoods of unspecified points. Points are derivable as perfect ideals of neighborhoods. The need for a neighborhood-based view of space arose in algebraic geometry and led to the geometric logic of Grothendieck toposes, showing that points, as unobservable particles that occupy no space but span the theory of space, are a convenient fiction. It took a long to notice this logical elephant that filled the room.

Another convenient fiction are the objects in category theory. The categorical view of mathematical structures focuses on morphisms and their compositions. The very first paper in category theory mentioned that the objects could be eliminated. Their role is to specify which morphisms are composable. They are the black boxes enclosing structures that are accessed along the morphisms as interfaces.  The enclosed structures may function identically along the isomorphisms, but we are expressly prohibited from making statements that involve objects' identities or equalities. Such statements have been called  \emph{evil}, with a tongue in cheek, but pointedly and objectively. The prohibitions of points and objects have a moral aspect. 

A holiday quest for a categorical approach to categories as black boxes led to \emph{catales}\/ as categories reduced to partial semigroups of morphisms,  analogous to locales as spaces reduced to lattices of neighborhoods. Just like points can be reconstructed as perfect filters of neighborhoods in locales, objects can be reconstructed as the top elements of the lattices of idempotents in catales. Interestingly, a precondition for the reductions and the reconstructions is that the categories involved are skeletal, i.e. that they satisfy the univalence axiom: that any isomorphic objects must be equal.
\end{abstract}

\newpage
\tableofcontents

\pagebreak

\section{Introduction: Geometry of black boxes}\label{Sec:intro}
% !TEX root = Z0-catale.tex

\subsection{Waves of arrows and particles of objects}
I recently collaborated on several papers with lots of idempotents \cite{PavlovicD:dede,PavlovicD:nuc}, and also wrote a book about categories where all objects arise by splitting the idempotents on a universal one \cite{PavlovicD:MonCom}. Explaining how particles arise by splitting waves, I came up with the picture of rocks and starlings at the end of this paper. It is probably not a very good metaphor but I tried to explain it here anyway. 

Category theory prides itself with studying mathematical structures on their own, abstracting away from irrelevant implementation details. E.g., a group may arise by permuting the zeros of a polynomial, or by transforming manifolds. The categorical view hides into black boxes the building materials of such groups and provides access to the group structure along the homomorphisms as explicit interfaces. A group boils down to how other groups see it, and how it sees other groups.\footnote{Category theory remains  perfectly capable of accounting for all concrete representations of any abstract group. See, e.g., \cite{CarboniA:matrices}. But it provides a handle on what is taken into account by distinguishing the abstract group structure from its concrete representations.} Yet the categories themselves are specified in terms of objects and morphisms, not the functors coming in and out. 

\subsection{Morality of categories}
One of the first principles of categories is that the equality between objects is \emph{evil}. Since objects are usually the first thing you define when you define a category, this is a bit like saying that one of the first principles of life is death. To such principles some people respond with humor, others with religion.

The reason why we cannot tell when two objects are equal is that they are black boxes: we cannot see inside.  We interact with them from outside, along morphisms as explicit  interfaces. Two black boxes may behave the same and look equal from outside but have different implementations inside. Integer arithmetic is implemented in many different ways in different computers. Many genetically different organisms (say, the grey wolf and the marsupial wolf) evolve and implement identical behaviors in different ways. Isomorphic shoes are often built differently in different factories. 

Category theory is a tool for studying concepts and structures \emph{on their own}, abstracting away from implementations. The task of abstracting away the implementations has permeated mathematics from the outset. Geometric constructions with ruler and compass have been taught by drawing particular triangles on paper, blackboard, in sand, or on a computer screen, but the requirement is that they should be valid \emph{on their own}, for all triangles on all media, irrespective of the implementations. The geometers are required to imagine that the lines drawn using a ruler are perfectly straight, infinitely long, infinitely thin, and that their intersection points are infinitely small and perfectly indivisible.

Why do we need abstractions? Chomsky says they are the essence of our language \cite{ChomskyN:cartesian}. They are necessary for communication. You say you love me. I say I love you. Do we mean the same thing? One website recommends yellow submarines. Another website also recommends yellow submarines. Do they sell the same thing? How can we know whether the same words refer to the same concept?

At the first sight, a concrete chair looks simpler than an abstract chair, and we try to agree about the abstract concept of a chair by pointing to a concrete chair. But when you take a closer look, the concrete chair is comprised of many parts and relationships, and it may not be clear which of the implementation details are relevant and which aren't. We may be able to learn what is a chair, love, a yellow submarine if we can interact with them. The observer observes them and accumulates some data. The scientist abstracts from the observations a hypothetical explanation and expresses it as an equation. The mathematician solves the equation. The solution explicates the information that was implicit in the equation, and the equation encodes the observation. But the explicit form of the solution allows computing predictions and testing the hypothetical explanation against observer's further observations. If needed, the scientist will abstract further explanations and the mathematician will solve further  equations. Concepts evolve through such cycles not only in science, but also in every baby's mind. Every baby contains an observer, a scientist, and a mathematician, of suitable age. Some of them may die as the baby ages.  A wide range of methods for concept learning has been implemented in recommender systems for the purposes of web commerce and social networking. An even wider range of methods for extracting explicit solutions from implicit models has evolved in science and engineering. Categories provide tools for structuring data and abstracting concepts. Categories provide frameworks for conceptualizing structures just like algebras provide frameworks for solving equations. Just like geometry as the science of space are simplified using points that occupy no space, but span it, category theory as the science of structure is simplified using objects that enclose structure in black boxes. But at the end of the day, the points must be reduced to the observable neighborhoods and the objects to waves of arrows.

\subsection{Idea}
Catales are related to categories like locales are to topological spaces:
\bear
\frac{\mbox{catales}}{\mbox{categories}} & \approx & \frac{\mbox{locales}}{\mbox{spaces}}
\eear
Intuitively, a locale is a lattice of open neighborhoods of a topological space, obtained by forgetting the points. Similarly, a catale is a partial semigroup of morphisms of a category, obtained by forgetting the objects.

\section{Locales and spaces}%: Pointless spaces}
\label{Sec:points}
% !TEX root = Z0-catale.tex

What are the advantages and what the disadvantages of the element-based approach of set theory, and of the morphism-based approach of category theory?

On one hand, the element-based approach has demonstrated its intuitive appeal and broad convenience ever since it was offered as the foundation of mathematics at the turn from XIX to XX century. On the other hand, with time, it seems to have acquired in mathematics the role that Latin used to play in philosophy: it provides a common ground for formal interactions in general councils, and an opportunity for highly abstract narratives, but when they go home, most people solve their problems using their native tongues.

The role of elements in mathematical reasoning is not unlike the role of points in geometry. On one hand, no one has ever seen a point, as it has no length, breadth, or depth, and it is by definition unobservable. When know that two intersecting lines determine a point, except that we shouldn't be able to see lines either, since they are infinitely thin, since they have no breadth or depth. We see surfaces, intersecting surfaces determine lines, and intersecting lines determine points. But we think in terms of points, because we think abstractly. And just like geometry studies space in terms of lines and points, which occupy no space and cannot be seen, set theory simplifies structures in terms of sets and elements, which have no structure. Simplification through abstraction. 

Every once in a while, though, abstract theory and concrete practice diverge, and need to be brought together. In algebraic geometry, the correspondence between the permutations of the zeros of a polynomial and the algebraic extensions of the ground field induced by them, observed by \'Evariste Galois, is extended to a correspondence between the sets of zeros of systems of polynomial equations and the sets of polynomials zeroed on given sets of tuples of elements. While the former sets are construed as \emph{algebraic}, since they are defined by polynomials, the latter are construed as \emph{geometric}, since they can can be thought of as topological neighborhood systems. Hence the vision of using the geometric abstraction of topology in algebra. \emph{Except}\/ that the obtained topological spaces, which came to be called \emph{spectral}\/ because of the analogy with the spectral methods used in studying the chemical elements --- turned out to often have geometrically indistinguishable points, i.e. such that every neighborhood of one point also contains the other point.

\subsection{Observations}
We call \emph{observations}\/ the logical interactions between a subject and an object. There is a wide gamut of situations that can be modeled in terms of observations. On one end of the gamut are observations as perceptions or measurements, where the observer may see, or hear, or touch an observable object and distinguish it from the environment. On the other end of the gamut are observations as tests or experiments, where the observer may empirically disprove a hypothesis or refute a conjecture about the object. Philosophical justifications for the prevalence of \emph{negative}\/ observations have been discussed at length in the theory of science \cite{FeyerabendP:against,LakatosI:refutations,PopperK:refutations}. Mathematical justifications are implicit in topos theory, and will be sketched below.

We present observations as elements of meet-semilattices.

\subsection{Equivalence of sober spaces and spatial locales}
Let $\Esp$ be the category of topological spaces and continuous maps between them. A topological space is a set $X$ with a topology, i.e. a family $\tp X \subseteq \WP X$ closed under finite intersections and all unions. The elements of a topology are thought of as open neighborhoods. A morphism $g\in \Esp(X,Y)$ is a function $g:X\to Y$ such that the induced inverse image $g^\ast:\WP Y \to \WP X$, defined by $g^\ast U = \{x\in X| gx \in U\}$, preserves the opens, i.e. maps any $U\in \tp Y$ into $g^\ast U\in \tp X$.

Let $\MLat$ be the category of meet-semilattices. An object $A$ is thus a partially ordered set with finite meets $\wedge\colon A\times A\to A$, capturing the greatest lower bounds. The top element $\top$ is construed as the empty meet. Given another meet-semilattice $B$, a morphism $f\in \MLat(A,B)$ is a function $f:A\to B$ which preserves the finite meets, i.e. $f(a\wedge c) = f(a) \wedge f(c)$ and $f\top = \top$. The elements of a meet-semilattice can be thought of as observations, or properties, and the meet $a \wedge c$ means that both properties $a$ and $c$ have been observed.

The set $2 = \{0,1\}$ lives in $\MLat$ as the ordering $\TTwo = \{0\lt 1\}$ and in $\Esp$ as the Sierpi\'{n}sky space $\Sier$ with the topology $\tp \Sier = \{1,2\}$. Since both $\Esp$ and $\MLat$ have products, we can define the powers 
\[ \Sier^A = \prod_{|A|} \Sier\   \mbox{ in } \Esp \qquad\qquad\mbox{and}\qquad\qquad \TTwo^X = \prod_{|X|} \TTwo\  \mbox{ in } \MLat\]  
where $|A|$ and $|X|$ denote the underlying sets. The inclusions $\MLat(A,\TTwo) \subseteq \Sier^A$ and $\Esp(X,\Sier) \subseteq \TTwo^X$
make $\pt A = \MLat(A,\TTwo)$ and $\tp X = \Esp(X,\Sier)$ into functors, as summarized on the following diagram. 
\beq\label{eq:galois-cat}
\begin{tikzar}[row sep = large]
X \ar[mapsto]{dd} \& \Esp \ar[bend right=15]{dd}[swap]{\tp}\ar[phantom]{dd}[description]\dashv \& %\radj \Phi V =
\displaystyle \MLat(A,\TTwo)\subseteq \Sier^{|A|} 
\\
\\
\displaystyle \Esp(X,\Sier)\subseteq \TTwo^{|X|}
\& \MLat^\op \ar[bend right=15]{uu}[swap]{\pt}
\& A \ar[mapsto]{uu}
\end{tikzar}
\eeq
The abuse of notation whereby $\tp X$ refers both to the topology on $X$ and to the meet-semilattice $\Esp(X,\Sier)\subseteq \TTwo^{|X|}$ is justified by the fact that any topology is a complete lattice, thus certainly also a meet-semilattice, and that that the lattices $\tp X$ and $\Esp(X,\Sier)$ are isomorphic, since every  $u\in \Esp(X,\Sier)$ corresponds a unique $U \in \tp X$ as $U = u^\ast 1$. The functor $\tp:\Esp \to \MLat$ can thus be viewed as forgetting the underlying set of the topological space $X$ and mapping it to the meet-semilattice $\tp X$ of its open neighborhoods. It is the \emph{open neighborhoods}\/ functor.

On the other hand, every $p\in \MLat(A,\TTwo)$ corresponds to a unique  $P\subset A$, defined by $P = p^\ast 0$. It is again easy to see that every such $P$ satisfies
\beq\label{eq:perfilt}
\top \not \in P  \quad \mbox{ and } \quad a\leq c \in P  \implies  a\in P\quad \mbox{ and }\quad
a \wedge c \in P  \implies  a\in P\, \vee\, c \in P
\eeq
The other way around, every $P\subset A$ satisfying \eqref{eq:perfilt} induces a unique $p\in \MLat(A,\TTwo)$ with $p a = 0 \iff a\in P$. The sets $P$ satisfying \eqref{eq:perfilt} are conveniently construed to be the \emph{abstract points}\/ of $A$. The functor $\pt:\MLat \to \Esp$ is then called the \emph{points}\/ functor. Any meet-semilattice $A$ is thus mapped to a topological space whose points are
\bear
\pt A & = & \{P\subset A\ \mbox{ satisfying } \eqref{eq:perfilt}\ \}
\eear 
The family open neighborhoods $\to\pt A$ is comprised of the unions of the basic open neighborhoods induced by the elements of the meet-semilattice $A$. The adjunction in \eqref{eq:galois-cat} makes this precise. Its unit and counit\footnote{The adjunction counit for the adjunction $\tp\dashv \pt$ is comprised of components in the form $\tp\pt A \tto\varepsilon A$. But they are in the dual category $\MLat^\op$, and the $\MLat$-components are in the form $\tp\pt A \oot\varepsilon A$.} are 
\begin{align*}
X & \tto\eta\ \pt \tp X & \tp\pt A & \oot\varepsilon\ A\\
x & \mapsto \hat x & \check{a} & \mapfrom  a
\end{align*}
%\begin{align*}
%X & \tto\eta\ \pt \tp X & A & \tto\varepsilon\ \tp \pt A\\
%x & \mapsto \hat x & a & \mapsto  \check{a}
%\end{align*}
where the operations $\hat{(-)}$ and $\check{(-)}$ are defined by
\beq\label{eq:hatcheck}
 U \in \hat x \iff x\not \in U\qquad \qquad\qquad P\in \check a \iff a\not \in P
 \eeq
Since for any $P\in \pt A$ by \eqref{eq:perfilt} holds
\[P\in \check{\overbrace{a\wedge b}} \ \iff\ a\wedge b \not\in P\ \iff\ a\not \in P \wedge b\not \in P\ \iff\ P \in \check{a} \cap \check{b}\]
the counit $\varepsilon$ maps $A$ into a basis of a topology, which is then generated by unions
\bea\label{eq:tppt}
\tp\pt A & = & \left\{ \bigcup_{a \in W} \check{a}\ |\ W\subseteq A\right\}
\eea
The isomorphisms $\tp X \stackrel{\eta}{\cong} \tp\pt \tp X$ and $\pt A \stackrel{\varepsilon}{\cong} \pt\tp\pt A$ confirm that both $\pt\tp:\Esp \to \Esp$ and $\tp\pt:\MLat\to \MLat$ are idempotent monads. Their algebras, \emph{viz}\/ fixpoints, form dual reflexive subcategories $\Esp^{\pt\tp}$ and $\MLat^{\tp\pt}$.
\beq\label{eq:galois-sober}
\begin{tikzar}[column sep=large,row sep = large]
\Esp \ar[bend right=15]{dd}[swap]{\tp}\ar[phantom]{dd}[description]\dashv 
\ar[bend right=10,two heads]{rr}\ar[phantom]{rr}[description]{\scriptstyle\top}
\&\& \Esp^{\pt\tp} 
\ar[bend right=15]{dd}[swap]{\tp}\ar[phantom]{dd}[description]\simeq 
\ar[bend right=10,hook]{ll}
\\
\\
\MLat^\op \ar[bend right=15]{uu}[swap]{\pt}
\ar[bend right=10,two heads]{rr}\ar[phantom]{rr}[description]{\scriptstyle\bot}
\&\& \left(\MLat^{\tp\pt}\right)^\op \ar[bend right=10,hook]{ll} \ar[bend right=15]{uu}[swap]{\pt}
\end{tikzar}
\eeq
The topological spaces contained in the category $\Ssp = \Esp^{\pt\tp}$ are called \emph{sober}. The meet-semilattices contained in the category $\MLat^{\tp\pt}$ are in fact complete, and satisfy 
\bear
a\wedge\big(\bigvee_{w\in W} w\big) &= & \bigvee_{w\in W} \big(a\wedge w\big)
\eear 
Viewed in  the opposite category $\MLat^\op$, such lattices are called \emph{frames}, or \emph{locales} \cite{JohnstoneP:stone,Joyal-Tierney:galois}.  The locales fixed by $\tp\pt$ are called \emph{spatial}, and they form the category $\Ssl = \left(\MLat^{\tp\pt}\right)^\op$. Since $\Ssp \simeq \Ssl$, the sober topological spaces can be studied as spatial locales, i.e. entirely in terms of their open neighborhoods --- in a point-free fashion \cite{JohnstoneP:pointless,JohnstoneP:elephant}. The  points are then derived as a second-order concept, as the sets of neighborhoods satisfying \eqref{eq:perfilt}.

\subsection{Non-examples}
\paragraph{A non-sober example.} Consider the topological space over the set $X= \NNn$ of natural numbers with the topology comprised of sets with a finite complement, i.e.
\bea\label{eq:tpNNn}
\tp \NNn & = & \left\{U\subseteq \NNn\ |\ \#\left(\NNn\setminus U\right) \lt \infty \right\}
\eea
Then any $n\in \NNn$ induces a point $\hat n\subset \tp \NNn$ defined in \eqref{eq:hatcheck}, i.e. $\hat n= \big\{P\subseteq \NNn\setminus \{n\}\big\}$. Hence the injection $\NNn \tto \eta \pt \tp \NNn$. For any $P\in \pt \tp \NNn$, the intersection $\bigcap P$ must have a finite complement. If this complement contains two elements, say $m$ and $n$, then $P$ would contain $\NNn\setminus \{m,n\}$ but not $\NNn\setminus \{n\}$ or $\NNn\setminus \{m\}$ and would not satisfy the second clause of \eqref{eq:hatcheck}. The complement of $\bigcap P$ for a point $P\in \pt \tp \NNn$ may therefore contain at most one number $n\in \NNn$. But the assignment $\NNn \tto \eta \pt \tp \NNn$, representing the abstract points $P\in \pt\tp\NNn$ in the form $P = \hat n$, is not surjective because
\bear
\varpi & = & \tp \NNn \setminus \{\NNn\}\ \ =\ \  \left\{U\subseteq \NNn\ |\ 0\lt \#\left(\NNn\setminus U\right) \lt \infty \right\}
\eear
is a point, i.e. $\varpi \in \pt\tp\NNn$, but $\hat n \neq \varphi$ for all $n\in \NNn$. So there is an abstract point in the space $\pt\tp\NNn$ that does not correspond to a concrete point of the space $\NNn$. However, $\varpi$ is an isolated point of $\pt\tp\NNn$, since none of the open neighborhoods in $\tp\pt\tp \NNn$, defined as in \eqref{eq:tppt}, or of the basis opens $\check U$, induced by $U\in \tp\NNn$ in \eqref{eq:tpNNn}, contains $\varpi$. Hence $\tp \NNn \stackrel{\tp\eta}\cong \tp\pt\tp \NNn$. This means that the continuous maps $\NNn\to \Sier$ have unique extensions to $\pt\tp\NNn\to \Sier$, and the difference between the two spaces is not observable by the continuous maps out.

\paragraph{A non-spatial example.} Consider the meet-semilattice $A = \TTwo^\NNn \big/ \sim$
of binary sequences modulo the equivalence relation
\bear
u\sim v & \iff & \#\big\{n\in \NNn\ |\ u_n \oplus v_n = 1 \big\} \lt \infty
\eear
where $\oplus$ is the \emph{exclusive or}. In words, the sequences $u,v:\NNn \to \TTwo$ are equivalent if they differ on at most finitely many indices $n\in \NNn$. Using $x\wedge(y\oplus z) = (x\wedge y)\oplus(x\wedge z)$, it is straightforward to show that $u\sim s$ and $v\sim t$ implies $u\wedge s \sim v\wedge t$, which means that
\bear
[u]_\sim \wedge [v]_\sim & = & [u\wedge v]_\sim
\eear
holds and makes $A$ into a meet-semilattice, where $[u]_\sim, [v]_\sim\in A$ denote the $\sim$-equivalence classes of $u$ and of $v$ respectively. The top and the bottom of $A$ are 
\bear
\top  & = & \Big\{u \in \TTwo^\NNn\ |\ \#\{n| u_n = 0\}\lt \infty\Big\}\\
\bot & = & \Big\{u \in \TTwo^\NNn\ |\ \#\{n| u_n = 1\}\lt \infty\Big\}
\eear
On the other hand, all $u:\NNn\to \TTwo$ with infinitely many 0s and 1s are nontrivially decomposable using $u = (u\vee v) \wedge (u\vee \neg v)$, and noticing that $s\lt t$ in $A$ means that $s$ has infinitely many 0s where $t$ has 1s. There are thus only two morphisms $A\to \TTwo$, one constantly $1$, and the other one taking only $\top$ to $1$ and everything else to $0$. The latter is $A \tto\varepsilon \tp \pt A$, since $\pt A = \Sier$. This time $\pt A \stackrel{\pt\varepsilon}\cong \pt\tp\pt A$ can be interpreted as saying that none of the observations $[u]_\sim \in A$ can be distinguished by morphisms $A\to \TTwo$, except $\top$ and $\bot$.

\subsection{Soberification and spatialization as absolute completions}
A limit or colimit in a category is called \emph{absolute}\/ if it is preserved by all functors. Bob Par\'e has proved in his 1969 thesis that the only absolute limits and colimits are the idempotent splittings.\footnote{Define the splitting as a triangle and specify the preorder of triangles?}

In the same vein, we say that 
\begin{itemize}
\item the spatial locale $\tp\pt X$ is the absolute completion of the locale $X$, 
\item the sober space $\pt\tp A$ is the absolute completion of the space $A$, \begin{itemize}
\item since the continuous maps out of $A$ uniquely determine the continuous maps out of $\pt\tp A$
\end{itemize}
and that
\begin{itemize}
\item since the localic maps into $X$ uniquely determine the localic maps into $\tp\pt X$,
\end{itemize}
\end{itemize}
If a space $X$ is not sober, its soberification $\pt\tp X$ has more points. However, the isomorphism $\tp X\stackrel{\tp \eta}\cong \tp\pt\tp X$ boils down to a bijection between the continuous maps $X\to \Sier$ and the continuous maps $\pt\tp X \to \Sier$. The space $\pt\tp X$ is thus the \emph{absolute completion}\/ of the space $X$, i.e. the completion by those points that must be preserved by any continuous map.

Similarly, if a meet-semilattice $A$ is not a spatial locale, the induced spatial locale $\to\pt A$ is a larger poset. However, the isomorphism $\pt A\stackrel{\pt\varepsilon}\cong \pt\tp\pt A$ means that any meet-preserving map $A\to \TTwo$ completely determines the extension $\tp\pt A \to \TTwo$. The locale $\tp\pt A$ is thus the \emph{absolute completion}\/ of $A$ because the elements added in $\tp\pt A$ must be preserved by any given morphism out of $A$.

Eliminating the evil unobservable points\footnote{with the same closures: Grothendieck!!} and reducing them to the points observable from their open neighborhoods is achieved through absolute completions, i.e. by adjoining the points that are observable by morphisms but not represented in the given space.

\section{Partial semigroups and categories}
 \label{Sec:psemi}
% !TEX root = Z0-catale.tex

\subsection{Categories}
A small category can be viewed as a diagram in the category $\Set$ of sets and functions
\bea\label{eq:Cat}
\CCc & = & \left(\begin{tikzar}[column sep = 5em]
C_0
\ar{r}[description]{\iota}
\&
C_1 \ar[shift left=1.5ex]{l}[description]{\varrho} \ar[shift right=1.5ex]{l}[description]{\delta}
\&
C_2
\ar{l}[description]{\mu}
\end{tikzar}\right)
\eea
where $C_0$ is the set of objects, or 0-cells, $C_1$ the set of morphisms, or 1-cells. The map $\iota$ assigns to each object its identity morphism, $\delta$ and $\varrho$ assign to each morphism its domain and codomain objects. The set $C_2$ of 2-cells is in the case of categories the pullback of $\delta$ and $\varrho$, i.e.
\bea
C_2 & = & \{<g,f>\in C_1\times C_1\ |\  \varrho (g)\oot{\ g\ } \delta (g)=\varrho (f) \oot{\ f\ } \delta (f)\}
\eea
The map $C_2\tto \mu C_1$ maps each $<g,f>\in C_2$ to the composite $\varrho( g) \oot{\ \mu(f,g)=g\circ f\ } \delta (f)$.

\paragraph{Objects as points.} Intuitively, the objects in $C_0$ can be thought of as  "points" of the category $\CCc$, and the morphisms in $C_1$ as the "open neighborhoods", or "observables". And just like some points of $X$ may be indistinguishable through observations, viewed as continuous functions $X\to \Sier$ (corresponding to the elements of $\tp X$), some objects of $\CCc$ may be indistinguishable through functors $\CCc\to \Set$.

\paragraph{Functors.} A functor $F\colon \CCc\to \DDd$ can be viewed as a diagram in $\Set$
\beq\label{eq:Fun}
\begin{tikzar}[column sep = 5em,row sep = 8ex]
C_0
\ar{r}[description]{\iota} \ar{d}[description]{F_{0}}
\&
C_1 \ar[shift left=1.5ex]{l}[description]{\varrho} \ar[shift right=1.5ex]{l}[description]{\delta} \ar{d}[description]{F_{1}}
\&
C_2 \ar[dashed]{d}[description]{F_{2}}
\ar{l}[description]{\mu}
\\
D_0
\ar{r}[description]{\iota} 
\&
D_1 \ar[shift left=1.5ex]{l}[description]{\varrho} \ar[shift right=1.5ex]{l}[description]{\delta} 
\&
D_2 
\ar{l}[description]{\mu}
\end{tikzar}
\eeq
where the squares formed by the arrows with the corresponding names commute and the function $F_{2}$ is induced by the \emph{functoriality condition}\/ on $F_{1}$:
\bear 
\varrho_{D}\circ F_{1}\circ \mu_{C} & = & \delta_{D}\circ F_{1}\circ \mu_{C}
\eear
and the fact that $D_{2}$ is a pullback of $\varrho_{D}$ and $\delta_{D}$. In terms of the elements $<g,f>\in C_{2}$, the functoriality condition is
\bear
F_{1}(g\circ f) & = & F_{1}(g) \circ F_{1}(f)
\eear

\subsection{Partial semigroups}
A \emph{semigroup}\/ structure on $A$ is an associative  binary operation $\ccirc$.
%, i.e. satisfying $(a \ccirc b)\ccirc c = a\ccirc (b\ccirc) c$. 
The associativity means that the first square in the following diagram commutes.
\beq\label{eq:PSg}
\begin{tikzar}{}
A\times A\times A \ar{rr}{A\times\ccirc}\ar{dd}[swap]{\ccirc\times A} \&\& A\times A \ar{dd}{\ccirc}
\&\&\&
A\times A\times A \ar{dd}[swap]{\ccirc \times A} \&\& A\times A \ar{ll}[swap]{A\times \widetilde\ccirc} \ar{dd}{\ccirc}
\\ \\
A\times A\ar{rr}[swap]{\ccirc} \&\& A 
\&\&\& 
A\times A  \&\& A  \ar{ll}{\widetilde\ccirc} 
\\ \\ \\
A\times A\times A \ar{rr}{A\times\ccirc} \&\& A\times A 
\&\&\&
A\times A \ar{rr}{\ccirc} \&\& A 
\\ \\
A\times A\ar{rr}[swap]{\ccirc} \ar{uu}{\widetilde\ccirc \times A} \&\& A \ar{uu}[swap]{\widetilde\ccirc} 
\&\&\&
A \ar{uu}{\widetilde \ccirc} \ar[equals]{uurr}[swap]{\id}
\end{tikzar}
\eeq
If the semigroup is considered in the category $\Rel$ of sets and relations, i.e. as a relation $\ccirc\in \Rel(A\times A, A)$, then there is also a dual relation $\widetilde \ccirc\in \Rel(A,A\times A)$. If $\ccirc$ is viewed as a subset of $(A\times A)\times A$ then $\widetilde \ccirc$ is the subset of $A\times (A\times A)$ defined by
\bear
<c,a,b>\in \widetilde\ccirc  & \iff & <a,b,c>\in \ccirc
\eear
Remembering that the relational composition of $R\in \Rel(X,Y)$ and $S\in \Rel(Y,Z)$ is the relation $S\circ R\in \Rel(X,Z)$ such that 
\bear
<x,z>\in S\circ R & \iff & \exists y.\ <x,y>\in R\  \wedge\  <y,z>\in S
\eear
it is easy to see that the commutative triangle on the  bottom-right in \eqref{eq:PSg}  says that for all $c,c'\in A$ \bear
\exists a, b \in A.\ <a,b,c> \in \ccirc \ \wedge\  <a,b,c'> \in \ccirc & \iff c = c' 
\eear
In other words, a relation $\ccirc\in \Rel(A\times A, A)$ that makes the triangle in \eqref{eq:PSg} commute is a surjective partial function $\ccirc:A\times A\pfn A$. We can thus write 
\[<a,b,c>\in \ccirc   \quad \mbox{ as } \quad a\ccirc b = c \qquad \mbox{ and }\qquad  
\exists c.\ <a, b, c> \in\ccirc \quad\mbox{ as } \quad a\ccirc b \halts 
\]
%
%
%
%
%When no confusion seems likely, we simplify notation even further and write
%\[
%a\ccirc b \quad\mbox{as} \quad ab
%\]
%In algebraic notation, the commutative square on the top-left of \eqref{eq:PSg} says that the partial function $\ccirc$ is associative, 
%\bea
%(a\ccirc b)\ccirc c & = &  a\ccirc (b\ccirc c)
%\eea whenever either side is defined. 
In this notation, the remaining two commutative squares in \eqref{eq:PSg},  on the top right and on the bottom left,  say
\bea \label{eq:frobcomp}
a\ccirc b\halts\  \wedge\ \  b\ccirc c\halts & \iff & (a\ccirc b)\ccirc c\halts\  \wedge\  \ a\ccirc (b\ccirc c)\halts
\eea   
These are the \emph{Frobenius}\/ conditions for the semigroup $A$ \cite{Carboni-Walters}. In ordered structures, they capture a version of modularity \cite{PavlovicD:QPL16}, or the presence of adjoints with respect to the order \cite{PavlovicD:spider}. The associativity square on the top left, of course, says that whenever the Frobenius conditions are met, the two composites must be equal, and the summary of \eqref{eq:PSg} is thus
\bea \label{eq:assoc}
a\ccirc b\halts\  \wedge\ \  b\ccirc c\halts & \iff & (a\ccirc b)\ccirc c\halts\  \wedge\  \ a\ccirc (b\ccirc c)\halts\ \wedge \ (a\ccirc b)\ccirc c  =  a\ccirc (b\ccirc c)
\eea 
The triangle at the bottom right, which makes the relation $\ccirc$ into a surjective partial function, is expressed in \eqref{eq:assoc} tacitly, by the use of the function notation. It makes the relation $\ccirc$ into a \emph{special}\/ Frobenius semigroup. Expanding the compositions in \eqref{eq:assoc} from binary and ternary to arbitrary finite lengths, as suggested by the associativity, one can show that the special Frobenius associativity in \eqref{eq:PSg} means that a composite is defined if and only if all of its decompositions are defined, i.e.
\begin{multline}
a_0\ccirc a_1\ccirc a_2\ccirc \cdots \ccirc a_n \halts \ \  \iff \\
 a_0\ccirc\cdots \ccirc a_{n_1}\halts \ \wedge\ \  a_{n_1}\ccirc\cdots \ccirc a_{n_2} \halts \ \wedge\ \  a_{n_2}\ccirc\cdots \ccirc a_{n_3} \halts \ \wedge\ \ \ldots \\
 \ldots \ \wedge\ \  a_{n_{m-2}}\ccirc\cdots \ccirc a_{n_{m-1}} \halts \ \wedge\ \  a_{n_{m-1}}\ccirc\cdots \ccirc a_{n_m}\halts
\end{multline}
where $0\lt n_1 \lt n_2 \lt \cdots \lt n_{m-1} \lt n_m = n$ is an arbitrary partition of the interval $[0,n]$. This is the \emph{spider form}\/ of the special Frobenius condition \cite{PavlovicD:spider}.

\subsection{The adjunction between categories and partial semigroups}\label{Sec:adj}
The adjoint functors are
\beq\label{eq:cat-catL}
\begin{tikzar}[row sep = large]
C_0\leftleftarrows C_1 \stackrel\circ \leftarrow C_2 \ar[mapsto]{dd}
 \& \Cat \ar[bend right=15]{dd}[swap]{\ootp}
 \ar[phantom]{dd}[description]\dashv \& 
A_0\leftleftarrows A_1 \leftarrow A_2 
\\
\\
\displaystyle C_1 \stackrel\ccirc \leftharpoonup C_1\times C_1
\& \PSg \ar[bend right=15]{uu}[swap]{\uupt}
\& \displaystyle A \leftharpoonup A\times A \ar[mapsto]{uu}
\end{tikzar}
\eeq
where
\begin{align}
&\qquad\qquad\qquad\qquad & A_0 & =  \{a\in A\ |\ a\ccirc a = a \}\notag\\
g\ccirc f & = \begin{cases}g\circ f & \mbox{ if }\delta_g = \varrho_{f}\\
\uparrow & \mbox{ otherwise}\end{cases}  \qquad& A_1 & =  \coprod_{a,b\in A_0}  A_{ab}\label{eq:idem-hom}\\
 && A_2 & =  \coprod_{a,b,c\in A_0} A_{ab} \times A_{bc}\notag
\end{align}
for $A_{ab} = \left\{f\in A\ |\  b\ccirc f\ccirc a = f\right\}$.
In other words,
\begin{itemize}
\item $\otp\CCc$ is the partial semigroup of morphisms of the category $\CCc$, and
\item $\upt A$ is the category whose objects are the idempotents of the partial semigroup $A$ and whose morphisms preserve the idempotents, in the sense that $f\in \pt A(a,b)$ satisfies $b\ccirc f = f = f\ccirc a$, or equivalently $b\ccirc f \ccirc a=f$.
\end{itemize}
The natural bijection
\bea\label{eq:adjunction}
\Cat\left(\CCc, \upt A\right) & \cong & \PSg\left(\otp \CCc, A\right)
\eea
is explained by the following diagram:
\[\begin{tikzar}[column sep = 2.5em]
C_0 \ar{rrr}[description]{\iota} \ar{ddd}[swap]{F_0} 
\&\&\& 
C_1 \ar[shift right=1.5ex]{lll}[description]{\delta} 
\ar[shift left=1.5ex]{lll}[description]{\varrho}\ar{ddd}[swap]{F_1} \ar{ddr}{F'}
\\ 
\\ 
\&\&\&\& A
\\
A_0 \ar{rrr}[description]{\iota} \&\&\& 
A_1\ar[shift right=1.5ex]{lll}[description]{\delta} 
\ar[shift left=1.5ex]{lll}[description]{\varrho} \ar[two heads]{ur}[swap]{\varepsilon}
\end{tikzar}
\]
A functor $F\in\Cat(\CCc,\upt A)$, with the object part $F_0$ and the arrow part $F_1$,  induces $F'\in \PSg(\otp \CCc, A)$ defined  
\bear
F' & = &  \varepsilon\circ F_1
\eear
where $\varepsilon$ projects $<a,b,f>$ from $A_1 = \coprod_{a,b\in A_0}  \left\{f\in A\ |\  b\ccirc f\ccirc a = f\right\}$ to $f$ in $A$. The other way around, a semigroup homomorphism  $F'\in \PSg(\otp \CCc, A)$ induces $F\in\Cat(\CCc,\upt A)$ comprised of
\bear
F_0(x) & = & F'\left(\iota_x\right)\\
F_1(h) & = & \left< F'\left(\iota_{\delta h}\right), F'\left(\iota_{\varrho h}\right), F'(h)\right>
\eear
The fact that the two mappings form a natural bijection is obvious in one direction and easily verified in the other.
The unit and the counit of the adjunction $\ootp\dashv \uupt:\PSg\to \Cat$, which are in the form
\beq
\begin{tikzar}{}%[column sep = large,row sep = large]
\CCc \ar{rr}{\eta} \&\& \upt{\otp{\CCc}} 
\&\&\& \otp {\upt {A}} \ar{rr}{\varepsilon}\&\& A 
\\ 
C_0 \ar{rr}{\eta_0}
\&\& \left(\upt{\otp{\CCc}}\right)_0
\&\&\&
\otp{\upt{A}} \ar{rr}{\varepsilon} \&\& A
\\
C_1 \ar[shift left]{u}\ar[shift right]{u}
\ar{rr}{\eta_1}
\&\&  \left(\upt{\otp{\CCc}}\right)_1  \ar[shift left]{u}\ar[shift right]{u}
\&\&\&
\otp{\upt{A}}\times \otp{\upt{A}} \ar[rightharpoonup]{u} \ar{rr}{\varepsilon\times \varepsilon} \&\& A\times A \ar[rightharpoonup]{u}
\\
C_2 \ar{u}
\ar[dashed]{rr} 
\&\&  \left(\upt{\otp{\CCc}}\right)_2 \ar{u}
\end{tikzar}
\eeq
are thus defined
\begin{align*}
\eta_0 x\ \  &=\ \  \iota_x  & \varepsilon<a,b,f> & = f\\
\eta_1 h\ \  &=\  \ \left< \iota_{\delta h}, \iota_{\varrho h}, h\right> 
\end{align*}
The category $\upt{\otp{\CCc}}$ is thus the absolute completion of $\CCc$. Each idempotent $a$ of $\CCc$ induces an object in $\upt{\otp{\CCc}}$, which plays the role of a splitting of $a$ viewed as a morphism. In this way, $\upt{\otp{\CCc}}$ splits the idempotents of $\CCc$, and turns out to be its absolute completion. On the other hand, the idempotents of the partial semigroup $A$ are turned into objects in the category $\upt A$, which get forgotten in the semigroup $\otp{\upt A}$, but remain recorded in the fact that every element $h$ of $A$ induces in $\otp{\upt A}$ as many triples $<a,b,h>$ as there are pairs of idempotents $a,b$ in $A$ such that $b\ccirc f\ccirc a = f$. The semigroup operation of $\otp{\upt A}$ is induced by the composition operation of the category $\otp A$, which means that $<a,b,h>$ and $<d,c,k>$ are composable if and only if $b=d$. The elements $<a,a,a>$ of $\otp{\upt A}$, induced by the idempotents $a$ from $A$, playing the role of identities in the category ${\upt A}$, can still be recognized in $\otp{\upt A}$ as the maximal elements with respect to the order of idempotents. This will allow us to characterize the partial semigroups that correspond to the morphism classes of categories.

\section{Taut categories and catales} \label{Sec:catale}
\subsection{Taut categories}
% !TEX root = Z0-catale.tex

\subsubsection{Absolute completion}
An endomorphism $\varphi\in \CCc(x,x)$ is said to be \emph{idempotent}\/ if $\varphi\circ \varphi = \varphi$. An \emph{idempotent splitting}\/ is usually defined as the decomposition
\[
\begin{tikzar}{}
\& y \ar{rr}[description]{id} \ar[tail]{dr}[description]{i}\&\& y\\
x \ar{rr}[description]{\varphi}\ar[two heads]{ur}[description]{q} \&\&x \ar[two heads]{ur}[description]{q} 
\end{tikzar}
\]

\begin{definition}\label{Def:CatA}
A category is \emph{absolutely complete} if all of its idempotents split, i.e.
\bear
\varphi\circ\varphi = \varphi & \implies & \exists i, q.\ \varphi = i\circ q\wedge q\circ i = \id
\eear
The category of absolutely complete categories is called $\CatA$.
\end{definition}

\begin{proposition}
The inclusion of absolutely complete categories into all categories has a left adjoint $\Karr\colon \Cat\to \CatA$, where $\Karr\CCc$, often written $\Kar\CCc$, has the idempotents of $\CCc$ as objects and the idempotent-preserving morphisms of $\CCc$ as the morphisms:
\bea
\lvert \Kar\CCc\rvert & = & \coprod_{x\in\lvert\CCc\rvert}\left\{\varphi_{x}\in \CCc(x,x)\ |\ \varphi\circ\varphi = \varphi\right\}\notag \\
\Kar\Ccc(\varphi_{x}, \psi_{y}) & = &\left\{f\in \CCc(x,y)\ |\ \psi_{y}\circ f \circ \varphi_{x} = f\right\} \eea
\end{proposition}

\paragraph{Idempotents are preordered.} The idempotents on $x$ obviously form a submonoid of the monoid $\CCc(x,x)$. The relation defined  on idempotents $\varphi, \psi$ on $x$ by
\bea
\varphi \leq \psi & \iff & \varphi = \varphi \circ \psi
\eea
is easily seen to be transitive and reflexive, i.e. a preorder. $\varphi\leq \psi$ and $\varphi\geq \psi$ imply $\varphi = \psi$ if and only if $\varphi \circ \psi = \psi \circ \varphi$. A convenient way to study the preorder of idempotents is to note that every idempotent forms the reflexive pair at the bottom of the following triangle.
\beq\label{eq:split-eq-coeq}
\begin{tikzar}[column sep  = large, row sep = 8ex]
\& y \ar[tail]{dr}[description]{i} \\
x \ar[two heads]{ur}[description]{q} \ar[shift right,bend right=30]{rr}[description]{\varphi} \&\& x \ar[shift left]{ll}[description]{id} \ar[shift right]{ll}[description]{\varphi} 
\end{tikzar}
\eeq

\begin{lemma} $i$ is the equalizer of the reflexive pair at the bottom of \eqref{eq:split-eq-coeq}, whereas $q$ is the coequalizer of the same pair.
\end{lemma}

\paragraph{Exercise.} If $\varphi\leq \psi$ then the splitting $y$ of $\varphi$ is a retract of the splitting $z$ of $\psi$,  which means that it can be obtained by splitting an idempotent on $z$.
\beq\label{eq:split-eq-coeq-preord}
\begin{tikzar}[column sep  = large, row sep = 8ex]
\&\& y \ar[tail]{ddrr}[description]{i} \ar[bend right,tail,dashed,thick]{d} \\
\&\& z \ar[tail]{dr}[description]{j} \ar[bend right,two heads,dashed,thick]{u} \\
x \ar[shift right,bend right=40]{rrrr}[description]{\varphi} \ar[two heads]{uurr}[description]{q} \ar{r}[description]{\varphi} \& x \ar[two heads]{ur}[description]{r} \ar[shift right,bend right=30]{rr}[description]{\psi} \&\& x \ar[shift left]{ll}[description]{id} \ar[shift right]{ll}[description]{\psi} \ar{r}[description]{\varphi} \& x
\end{tikzar}
\eeq

\subsubsection{Skeleton}
\begin{definition}\label{Def:CatS}
A category is \emph{skeletal} if all of its isomorphisms are automorphisms, i.e. any isomorphic object are equal:
\bea\label{eq:univalence}
A\cong B & \implies & A=B
\eea
The category of skeletal categories is called $\CatS$.
\end{definition}

\begin{proposition}
Assuming the axiom of choice, the inclusion of skeletal categories into all categories has a left adjoint  $\Skell\colon \Cat\to \CatS$, where $\Skell \CCc$, often written $\Skel\CCc$, consists of the isomorphism classes of $\CCc$-objects as $\Skel\CCc$-objects and arbitrary $\CCc$-morphisms between arbitrary elements of the isomorphism classes as $\Skel\CCc$-morphisms:
\bea
\lvert \Skel\CCc\rvert & = & \left\{\XXx\subseteq \lvert\widetilde \CCc\rvert\ |\ \forall x x' \in \XXx \exists \kappa \in \widetilde \CCc(x,x')\right\}\notag\\
\Skel\CCc(\XXx, \YYy) & = &\left\{f\in \CCc(x,y)\ |\ x\in \XXx \wedge y\in \YYy\right\}\ \Big/ \ \iapprox
\label{eq:Skel} 
\eea
where $\widetilde\CCc\inclusion \CCc$ is  the subgroupoid of isomorphisms 
\bea
\lvert \widetilde\CCc\rvert & = & \lvert \CCc\rvert\notag\\
\widetilde\CCc(x,y) & = &\left\{f\in \CCc(x,y)\ |\ \exists f'\in \CCc(y,x).\ f'\circ f = \id_{x}\wedge f\circ f' = \id_{y}\right\} 
\eea
and $\wwidetilde \CCc\inclusion \widetilde \CCc$ is an indiscrete subpreorder, i.e. a subcategory with a single morphism in each homset, which we denote $\iota_{xy}\in \wwidetilde\CCc(x,y)$ and define the equivalence relation in \eqref{eq:Skel}  by
\bear
\left(x\tto f y\right) & \iapprox & \left(x'\tto{\iota} x\tto f y \tto \iota y'\right)
\eear
\end{proposition}

\begin{lemma}\label{Lemma:groups}
For any connected component  $\XXx$ of $\widetilde \CCc$ there is a group $G_{\XXx}$ with
\begin{enumerate}[a)]
\item for every $x\in \XXx$ a canonical group isomorphism $\widetilde \CCc(x,x) \cong G_{\XXx}$, and
\item for all $x,x' \in \XXx$  and every $\kappa \in \widetilde \CCc(x,x')$ a unique bijection $\widetilde \CCc(x,x') \stackrel\kappa\cong G_{\XXx}$.
\end{enumerate}
\end{lemma}

\bpr \begin{enumerate}[a)]
\item Take $G_{\XXx}$ to be the group of automorphisms $\widetilde \CCc(x,x)$ for any $x\in \XXx$. For any $x'\in \XXx$, any pair $\kappa\in \widetilde \CCc(x,x')$ and $\kappa'\in \widetilde \CCc(x',x)$ where $\kappa'\circ\kappa = \id_{x}$ and $\kappa\circ\kappa' = \id_{x'}$ clearly establish a group isomorphism of $\widetilde \CCc(x,x)$ and $\widetilde \CCc(x',x')$.  All monoids $\widetilde \CCc(x',x')$ are thus isomorphic to $G_{\XXx}$.

\item For all $x,x'\in \XXx$ and $\kappa, \chi \in \widetilde \CCc(x,x')$ the isomorphism $[\chi,\kappa]  = \chi^{-1}\circ \kappa \in \widetilde \CCc(x,x)$ gives $\kappa = \chi\circ [\chi,\kappa]$. Fixing an arbitrary isomorphism $\kappa \in \CCc(x,x')$ fixes the bijection
\bea
\widetilde\CCc(x,x') & \tto{\kappa}  &\widetilde\CCc(x,x)\notag\\
\chi & \mapsto & [\chi,\kappa]\label{eq:kappa}
\eea
Postcomposing with $\CCc(x,x)\cong G_{\XXx}$ established in (a) yields the claimed bijection $\widetilde\CCc(x,x') \stackrel\kappa \cong G_{\XXx}$.
\end{enumerate}
\epr

\begin{corollary}
For every component $\XXx \in \lvert \Skel\Ccc\rvert$, the hom-set $\Skel\Ccc(\XXx, \XXx)$ is the group $G_{\XXx}$ from Lemma~\ref{Lemma:groups}. \end{corollary}

\begin{proposition}\label{Lemma:components}
Every choice of\  \ $\wwidetilde \CCc \inclusion \CCc$ induces an equivalence of categories $\varepsilon\colon \CCc \tto{\ \sim\ } \Skel \CCc$, mapping each object $x$ of $\CCc$ to its isomorphism class $\XXx$ and every morphism $f$ to its equivalence class $[f]_{\approx}$. 
\end{proposition}

\begin{proposition}
The skeleton of an absolutely complete category is absolutely complete.
\end{proposition}

\subsubsection{Taut completion}
\begin{definition}\label{Def:CaT}
A category that is both skeletal and absolutely complete is said to be \emph{taut}. The category of taut categories is called $\CaT$.
\end{definition}

\begin{proposition}
The inclusions of taut categories into all categories has a left adjoint $\Tautt\colon \Cat \to \CaT$, where $\Tautt \CCc$ is often written $\Taut\CCc$ is defined
\bear
\lvert \Taut\CCc\rvert & = & \lvert \Kar \CCc\rvert\  \big/  \approx \notag \\
\Taut\Ccc\left([\varphi], [\psi]\right) & = & \Kar \CCc(\varphi, \psi) \label{eq:Taut} 
\eear
where the equivalence relation $\approx$ characterizes the idempotents that are isomorphic in $\Kar\CCc$ by
\bea\label{eq:approx}
a\approx b & \iff & \exists uv.\ a = u\circ v\ \ \wedge\  \ b = v\circ u\ \ \wedge\  \ u\circ v\circ u = u\ \ \wedge\  \ v\circ u\circ v = v
\eea
and \eqref{eq:Taut} is sound by Lemma~\ref{Lemma:idem-iso}.
\end{proposition}

\paragraph{Condition \eqref{eq:approx} characterizes isomorphic idempotents.} The assumption $u\ccirc v\ccirc u = u$ implies that $a = u\ccirc v$ is idempotent, and $v\ccirc u\ccirc v = v$ implies the same for $b = v\ccirc u$.  Moreover, $a$ and $b$ are among the units of $u$ and $v$, since $u = a\ccirc u \ccirc b$ and $v = b\ccirc v \ccirc a$.  The other way around, $u$ and $v$ are morphisms between $a$ and $b$, i.e.
\beq\label{eq:isocatale}
\begin{tikzar}{}
a \ar[bend left]{rr}[description]{v} \&\& b\ar[bend left]{ll}[description]{u}
\end{tikzar}
\eeq
Last but not least, $u$ and $v$ make $a$ and $b$ isomorphic as idempotents, by the very definitions  $a = u\ccirc v$ and $b = v\ccirc u$. Any pair of idempotents $a$ and $b$ decomposed in the form $a=u\ccirc v$ and $b = v\ccirc u$ is thus isomorphic. The catale condition thus makes such isomorphic idempotents equal.

\begin{lemma}\label{Lemma:idem-iso}
If two idempotents $a = \left(A\eepi{v_{A}} X \mmono{u_{A}} A\right)$ and $b=\left(B\eepi{u_{B}} X \mmono{v_{B}} B\right)$ split on the same object $X$ in a category $\CCc$, then $u=u_{A}\circ u_{B}$ and $v = v_{B}\circ v_{A}$ make them isomorphic in $\Kar\CCc$.
\end{lemma}

\bpr
\beq\label{eq:iso-split}
\begin{tikzar}[column sep=2.5cm]
A  \arrow[loop, out = 140, in=220, looseness = 6]{}[swap]{a} 
\arrow[two heads,shift right = .5ex,shorten=1mm,bend right = 9]{r}[swap]{v_{A}}
\arrow[shift right = .5ex,shorten=1mm,bend right = 30]{rr}[swap]{v}    
\& 
X \arrow[tail,shift right = .5ex,shorten=1mm,bend right = 9]{r}[swap]{v_{B}} 
\arrow[tail,shift right = .5ex,shorten=1mm,bend right = 9]{l}[swap]{u_{A}}
\&
B
\arrow[loop, out = 40, in=-40, looseness = 6]{}{b} 
\arrow[shift right = .5ex,shorten=1mm,bend right = 30]{ll}[swap]{u} 
\arrow[two heads,shift right = .5ex,shorten=1mm,bend right = 9]{l}[swap]{u_{B}}
\end{tikzar}
\eeq
Since $u = a\circ u \circ b$ means $u\in \Kar\PPP(a,b)$, and $v = b\circ v \circ a$ justifies $v\in \Kar\PPP(b,a)$ and the equations $u\circ v = a$ and $v\circ u = b$ make them into an isomorphism, since the $\CCc$-idempotent $a$ is the identity on the $\Kar\CCc$-object $a$, whereas $b$ is the identity on $b$. 
\epr

\paragraph{Summary.} The relations between the categories of categories specified in Definitions~\ref{Def:CatA}, \ref{Def:CatS}, and \ref{Def:CaT} is
\beq\label{eq:Relations}
\begin{tikzar}[column sep=7ex,row sep = 8ex]
\& \Cat \ar[bend right=15,two heads]{dl}[swap]{\Karr}
\ar[bend right=15,two heads]{dr}[swap]{\Skell}
\ar[phantom]{dl}[description]\dashv
\ar[phantom]{dr}[description]\dashv 
\\
\CatA \ar[hook,bend right=15]{ur} \ar[two heads,bend right=15]{dr}[swap]{\Skell}
\ar[phantom]{dr}[description] {\dashv}
\&\& \CatS \ar[hook,bend right=15]{ul} 
\ar[two heads,bend right=15]{dl}[swap]{\Tautt}
\ar[phantom]{dl}[description] {\dashv} 
\\
\&
\CaT \ar[hook,bend right=15]{ur}\ar[hook,bend right=15]{ul}
\end{tikzar}
\eeq

\subsection{Catales}
% !TEX root = Z0-catale.tex

\begin{definition}
An element $a$ of a partial semigroup $A \stackrel\ccirc\leftharpoonup A\times A$ is called an\/ \emph{identity} if all $f,g\in A$ satisfy
\bea
f\ccirc a \halts & \implies &  f\ccirc a = f\\
a\ccirc g \halts & \implies & a\ccirc g = g
\eea
and moreover $a\ccirc f\halts$.
\end{definition}

\paragraph{Notation.} We write
\begin{itemize}
\item $\idds A$ to denote the set of all identities and
\item $\kar A$ to denote the set of all idempotents 
\end{itemize}
contained in the partial semigroup $A$.

\begin{lemma}\label{Lemma:idds} All identities $a,b$ satisfy the following implications:
\begin{enumerate}[a.]
\item $a\ccirc b\halts\ \implies\ a=b$
\item $a\ccirc f\halts \wedge\ \ b\ccirc f \halts \ \implies\ a=b$
\item $f\ccirc a\halts \wedge\ \ f\ccirc b \halts\  \implies\  a=b$
\end{enumerate}
\end{lemma}
\bpr
\begin{enumerate}[a.]
\item Since $a$ is an identity, $a\ccirc b\halts$ implies $a\ccirc b = b$. Since $b$ is an identity, $a\ccirc b\halts$ implies $a\ccirc b = a$. Hence $a = a\ccirc b = b$. 
\item Since $a\ccirc f\halts$ implies $a\ccirc f=f$, $b\ccirc f\halts$ implies $b\ccirc f = b\ccirc a\ccirc f\halts$. \eqref{eq:frobcomp} then implies $b\ccirc a\halts$ and (a) gives $a=b$.
\item is proved dually to (b).
\end{enumerate}
\epr

\begin{definition}\label{Def:catale}
A partial semigroup $A$ is called a \emph{catale} if it satisfies the following conditions
\begin{enumerate}[a)]
\item  for every $f\in A$ there are identities $a,b\in \idds A$ such that \begin{itemize}
\item $b\ccirc f\ccirc a\halts$
\end{itemize}
\item for every idempotent $\varphi\in \kar A$ there are $i,q\in A$ such that 
\begin{itemize}
\item $i\ccirc q = \varphi$, 
\item $q\ccirc i \in \idds A$,
\item if $j\ccirc r = i\ccirc q$ and $r\ccirc j \in \idds A$ then $r\ccirc j = q\ccirc i$.
\end{itemize}
\end{enumerate}
\end{definition}

\paragraph{Consequences of \textit{(a)}.} Def.~\ref{Def:catale}{\it (a)}\/ implies $b\ccirc f\halts$ and $f\ccirc a\halts$ by \eqref{eq:assoc}. The assumption that $a$ and $b$ are identities implies that $b\ccirc f=f = f\ccirc a$. Lemma~\ref{Lemma:idds} implies every $f$ in a catale determines a unique right identity $a$ and a uniqe left identity $b$. In the category of sets and functions, there is thus the diagram \bea\label{eq:IddsA}
\Idds A & = & \left(\begin{tikzar}[column sep = 5em]
\idds A \ar[tail]{r}[description]{\id}
\&
A \ar[shift left=1.5ex, two heads]{l}[description]{\varrho} \ar[shift right=1.5ex,two heads]{l}[description]{\delta}
\end{tikzar}\right)
\eea
where $\id$ is the inclusion. Can we extend it to \eqref{eq:Cat} and get a category? 

\paragraph{Every catale induces a category.} The identity $\delta_{f}$ satisfying $f\ccirc \delta_{f}=f$ is called the \emph{domain}\/ of $f$. The identity $\varrho_{f}$ satisfying $\varrho_{f}\ccirc f=f$ is called the \emph{codomain}\/ of $f$. Lemma~\ref{Lemma:idds} now implies $$\varrho\circ\id\  =\  \mathsf{id}_{\idds A}\  =\  \delta\circ \id$$ 
The decomposition
\bea\label{eq:ids-hom}
A & = & \coprod_{a,b \in \idds A} A_{ab} \mbox{ where}\\
&& A_{ab} = \{f\in A\ |\ \delta_{f}=a\wedge \varrho_{f} = b\}
\eea
makes every catale into the family of morphisms of a category $\Idds A$, as in \eqref{eq:IddsA}, over the family of objects  $\idds A$ as and the hom-sets $A_{ab}$.  \eqref{eq:IddsA} can be extended to \eqref{eq:Cat} using
\bea\label{eq:ids-hom-two}
A_{2} & = & \coprod_{a,b,c \in \idds A} A_{ab}\times A_{bc} 
\eea
The composition $A\oot\mu A_{2}$ is induced by the partial semigroup operation on $A$. The fact that it completely determines it easily follows from Def.~\ref{Def:catale}{\it (a)}. The following lemma is equally easy to prove, but we state it for the record.

\begin{lemma}\label{lemma:identities}
An element of a catale $A$ is an identity if and only if it is idempotent and equal to its domain and codomain
\bea
f\in \idds A & \iff& f = f\ccirc f = \delta_{f}= \varrho_{f}
\eea
\end{lemma}

\paragraph{Consequences of \textit{(b)}.} Def.~\ref{Def:catale}{\it (b)}\/ says that every idempotent has a splitting and that two identities that arise by splitting the same idempotent must be equal. The former claim means that the category $\Idds A$ is absolutely complete. The latter claim means that the category $\Idds A$ is skeletal. In summary, 

\begin{proposition}\label{prop:Idds-taut}
For any catales $A$, the category 
$\Idds A$ is taut. \end{proposition}

\begin{definition}
The category $\CaL$ of catale is the full subcategory of the category $\PSg$ of partial semigroups spanned by catales.\end{definition}

\paragraph{Remark.} Catale morphisms are thus the semigroup morphisms between catales. Note that they have to map idempotents to idempotents but not necessarily identities to identities. The idempotent to which an identity is mapped may not be an identity.

\subsection{The equivalence of taut categories and catales}\label{Sec:catcaleq}
% !TEX root = Z0-catale.tex

In Sec.~\ref{Sec:adj} we described the adjunction between categories and partial semigroups, displayed on the left in the next diagram.
\beq\label{eq:taut-catale}
\begin{tikzar}[column sep=large,row sep = large]
\Cat \ar{dd}[swap]{\ootp}\ar[phantom,bend left=20]{dd}[description]{\dashv} 
\&\& \CaT 
\ar{dd}[description]{\hspace{.75ex}\ootp}
\ar[phantom,bend right=22.5]{dd}[description]{\dashv} \ar[phantom,bend left=27.5]{dd}[description]{\dashv} 
\ar[hook']{ll}
\\
\\ 
\PSg \ar[bend right=37.5]{uu}[swap]{\uupt}
\&\& \CaL \ar[hook']{ll} \ar[bend left=35]{uu}{\IIdds}\ar[bend right=40]{uu}[swap]{\IIdds}
\end{tikzar}
\eeq
On the right, the category of categories  $\Cat$ is restricted to its full subcategory of taut categories $\CaT$, whereas the category of partial semigroups is restricted to its full subcategory of catales $\CaL$. It is straightforward to check that the functor $\ootp\colon \Cat\to \PSg$ restricts to $\ootp\colon \CaT\to \CaL$.

On the other hand, the functor $\uupt\colon \PSg\to \Cat$ restricts to $\IIdds\colon \CaL\to \CaT$ because all inclusions $\Idds A \hookrightarrow \upt A$ are equivalences, since by Prop.~\ref{prop:Idds-taut},  $\Idds A$ is absolutely complete. The adjunction $\ootp \dashv \uupt\colon \PSg\to \Cat$ from Sec.~\ref{Sec:adj} therefore restricts to $\ootp\dashv\IIdds \colon \CaL\to \CaT$. To see that this is an equivalence, we first spell out the natural bijection 
\bea\label{eq:adjunction-back}
\CaT\left(\Idds A, \CCc\right) & \cong & \CaL\left(A, \otp \CCc\right)
\eea
extablishing the inner adjunction $\IIdds\dashv \ootp\colon\CaT \to \CaL$ in \eqref{eq:taut-catale} on the right. The correspondence in \eqref{eq:adjunction-back} is  explained using the following diagram:
\beq\label{eq:expl}
\begin{split}
\begin{tikzar}[column sep = 2.5em]
A_0 \ar{ddd}[swap,pos=0.66]{G'_0} \ar{rrr}[description]{\id}
\&\&\& 
A_1 \ar[shift right=1.5ex]{lll}[description]{\delta} 
\ar[shift left=1.5ex]{lll}[description]{\varrho}\ar{ddd}[swap,pos=0.66]{G'_1} 
\\ 
\&\idds A \ar[crossing over]{rrr}[description]{\id} \ar[hook',shift left=1]{ul}\&\&\& 
A \ar[shift right=1.5ex,crossing over]{lll}[description]{\delta} 
\ar[shift left=1.5ex,crossing over]{lll}[description]{\varrho} \ar[shift right=1]{ul}[swap]{H}
\ar[shift left=0.75]{ddl}{G}
\\
\\
C_0 \ar{rrr}[description]{\id} \&\&\& 
C_1\ar[shift right=1.5ex]{lll}[description]{\delta} 
\ar[shift left=1.5ex]{lll}[description]{\varrho}
\end{tikzar}
\end{split}
\eeq
Like before, $A_0$ are the idempotents of $A$ and $A_{1}$ is the set of triples $<a,b,f>$ where $a$ and $b$ are idempotent and $b\ccirc f\ccirc a = f$. Set $Hf = <\delta_{f},\varrho_{f}, f>$.  Any $G\in \CaL\left(A, \otp \CCc\right)$ surely maps the elements of $\idds A$ to idempotents in $C_{1}$, but it is not obvious that it maps the identities of the catale $A$ to identities in $C_{1}$. It induces  a unique functor $G'\in \CaT \left(\upt A, \CCc\right)$ comprised of
\begin{itemize}
\item the object part $G'_{0}$, which maps $a\in A_{0}$ to the object that splits the idempotent $Ga$ in $\CCc$;
\item the arrow part $G'_{1}$ which maps $<a,b,f>\in A_{1}$ to the morphism from the splitting of $Ga$ to the splitting of $Gb$ induced by $Gf$. 
\end{itemize}
$G'_{0}$ and $G'_{1}$ are well-defined because $\CCc$ is taut, which means that every idempotent has a unique splitting. The fact that $G'_{0}$ and $G'_{1}$ form a functor follows directly from their definitions. The fact that $G= G'_{1}H$ preserves the identities follows from Lemma~\ref{lemma:identities}. We have thus established the bijective correspondence in \eqref{eq:adjunction-back}, leading to

\begin{theorem}\label{Thm:equivalence}
The functors $\begin{tikzar}[column sep = 5em]
\CaT \ar[bend right = 10]{r}[swap]{\ootp} \ar[phantom]{r}{\simeq} \& \CaL \ar[bend right = 10]{l}[swap]{\IIdds}
\end{tikzar}$ form an equivalence.
\end{theorem}

\bpr
The following diagram provides an overview of the discussion preceding the theorem.  
\beq\label{eq:cat-catL-back}
\begin{tikzar}[row sep = large]
\idds A\leftleftarrows A \leftarrow A_2 
  \& \CaT \ar[bend left=15]{dd}{\ootp}\ar[phantom]{dd}[description]\dashv \& 
C_0\leftleftarrows C_1 \leftarrow C_2 \ar[mapsto]{dd}
\\
\\
A \leftharpoonup A\times A
\ar[mapsto]{uu}
\& \CaL \ar[bend left=15]{uu}{\IIdds}
\&  \displaystyle C_{1} \stackrel\ccirc \leftharpoonup C_{1}\times C_{1}
\end{tikzar}
\eeq
The fact that this adjunction is an equivalence follows by inspection of the unit $A\stackrel\eta\cong \otp{\left(\Idds{A}\right)}$ and the counit $\Idds{\left(\otp{\CCc}\right)}\stackrel\varepsilon\cong \CCc$.
\epr

\section{Examples} \label{Sec:xample}
% !TEX root = Z0-catale.tex

\subsection{Folding and unfolding extensional types}
In set theory, sets are sets of elements and the sets $A$ and $B$ are equal if their elements are equal, i.e.  $A=B \iff  \forall x. x\in A\Leftrightarrow x\in B$. The same elements form the same set.  In category theory, objects constrain the inputs and the outputs of morphisms. In the category $\Set$, \emph{\textbf{sets  are reduced to types}}. The sets $A$ and $B$ are reduced to the input and the output types of the functions in the form $f\colon A\to B$, forming the hom-set $\Set(A,B)$. Since a set as an object is a black box, its elements can only be accessed through functions. Since we cannot tell what they consist of, we don't know if $x\in A$ is the same as $y\in B$, and can never tell if $A$ is the same as $B$. A categorical theory of sets is a theory of types.

The problem of assigning types to the terms of a type-free $\lambda$-calculus was considered by Dana Scott in \cite{ScottD:relating}. The idea was that the domain $U$ of a type-free $\lambda$-calculus should be viewed as a universal type giving rise to all other types is its retracts. The retraction of $U$ to a type $A$ thus gives rise to an idempotent $\alpha = \left(U\epi  A \mono  U\right)$, splitting on $A$. An untyped $\lambda$-term $f$ can be typed in the form $f\colon A\to B$ if $\beta \circ f\circ \alpha = f$ where $\beta = \left(U\epi B \mono U\right)$. This requirement determines the hom-sets of the category $\upt A$ induced by the partial semigroup $A$ in Sec.~\ref{Sec:adj}. The main technical problem addressed in \cite{ScottD:relating} was to encode the structure of cartesian closed categories, tightly corresponding to typed $\lambda$-calculi, in the monoid over the universal type $U$. Leaving that issue aside, any monoid or even a partial semigroup induces a category where its idempotents are the objects. The functor $\uupt\colon \PSg \to \Cat$ in Sec.~\ref{Sec:adj} summarizes this construction. The category $\set_{\leq \upsilon}$ of sets below some cardinal $\upsilon$ can be obtained by splitting the idempotents in the monoid $M_{\upsilon}= \Set(\upsilon,\upsilon)$ of functions on $\upsilon$. For every pair of retracts $A, B$ of $\upsilon$ every function $f$ on $\upsilon$ induces a unique function $f_{AB}$.
\[
\begin{tikzar}[sep=large]
\upsilon \ar{r}[description]{f} \ar[two heads,bend right=15]{d} \& \upsilon \ar[two heads,bend right=15]{d}\\
A \ar{r}[description]{f_{AB}} \ar[hook,bend right=15]{u} \& B \ar[hook,bend right=15]{u}
\end{tikzar}
\]
The objects of $\set_{\leq\upsilon}$ are retracts of $\upsilon$ and its morphisms are retracts of functions on $\upsilon$. As each idempotent on $\upsilon$ is taken as an identity, i.e. an object, all idempotents below it are copied as idempotents on that object. The elements of the monoid on $\upsilon$ are thus branched into many copies to form the category $\set_{\leq\upsilon}$. A universe $\Set$ of all sets can be derived in a similar way from the monoid of functions on an inaccessible cardinal. The monoid is then a proper class and does not live in the universe $\Set$. Nevertheless, the categories $\Set$ and $\set_{\leq\upsilon}$ can be folded back into the monoids from which they were unfolded, just like the typed $\lambda$-calculus and it cartesian closed category, constructed by Dana Scott from the monoid on the universal domain $U$, can be folded back to the untyped $\lambda$-calculus on $U$, with no loss of information.

Note that this is \emph{not}\/ how the equivalence $\CaT\simeq \CaL$ works in Sec.~\ref{Sec:catcaleq}. The elements of a catale are required to remember their identities. This may be elided by local encodings in special cases, but in general, the global folding of local information, like in categories of sets and in cartesian closed categories with a universal comain --- is \emph{not}\/ possible. It is only possible in \emph{extensional}\/ universes. To see this, take a Turing complete programming language $\PPp$. It comes with some native type system, usually some basic datatypes, their products and sums, polymorphic list constructors, and so on. Invariably, many more types can be defined by programming suitable type checkers. A type checker for type $A$ can be programmed as an idempotent function, filtering the inputs of type $A$ and delivering them at the output unchanged, and filtering out the inputs that are not of type $A$. For any pair of types $A,B$, the language $\PPp$ comes with a program evaluator that runs programs and computes functions of type $A\to B$. Programming and computation in terms of such abstract evaluators is described in some detail in \cite{PavlovicD:MonCom}. A computer is presented as a universe of types: a monoidal category. The fact that there many different programs evaluate on each computation means that the universe of types is \emph{intensional}\/ in general. The consequence is that the program interpreters for computations of type $A\to B$, for different $A$s and $B$s, are \emph{not}\/ retracts of the program interpreters of type $\PPp\to \PPp$, as they are in $\lambda$-calculus and in set theory, although all types $A$ and $B$ are retracts of $\PPp$. It was worked out at the end of chapters 2 and 8 that the interpreters of all types are retracts of the one of type $\PPp$ if and only if the type universe is essentially extensional, which means that it contains an extensional retract, where every computation has a unique program. That is why a catale of computations has to carry explicit type annotations as identities. Intuitively, this refers to the fact that program  interpreters, just like natural language interpreters, genuinely use type information, and interpret differently typed programs and statements differently.

\subsection{The catale of spaces is a space}
In a suitably large universe, the category $\Esp$ is small and the functor $\tp$ can be applied to it. The underlying set of the catale $\tp\Esp$ is the class $\Esp^\tp$ of all continuous maps $f:X\to Y$, for arbitrary topological speces $X$ and $Y$, modulo the congruence generated by identifying the idempotents $a:X\to X$ and $a':Y\to Y$ such that the retracts $A =\{x\in X\ |\ a(x) = x\}$ of $X$ and $A' = \{y\in Y\ |\ a'(y) = y\}$ of $Y$ are homeomorphic.

The catale $\tp\Esp$ thus forgets the identities of individual topological spaces and retains the continuous maps in its underlying class and their composability in its semigroup operation. The individual topological spaces are recorded as idempotents. Intuitively, the topological spaces that contain homeomorphic retracts are glued along the shared retract. Together, they can be thought of as a large topological space. This echoes Grothendieck's idea of the category of topological spaces as a generalized topological space. The idea of gluing a pair of topological spaces along their homeomorphic retracts echoes the concept of a sheaf. The category of topological spaces, viewed as a catale, can thus be construed as a large sheaf. The base space is implicit but can be reconstructed from path liftings. 

\subsection{The category of locales is the catale of neighborhoods}
Local homeomorphisms as etale spaces over any given space form a Grothendieck topos. Grothendieck's idea was, of course, that this topos is the generalization of the space, since the etales over a space externalize its geometry. Gluing together all toposes as fibers of a fibration over the category of all spaces yields what Grothendieck called the \emph{gros topos} \cite{GrothendieckA:SGA4,JohnstoneP:elephant}. Remarkably, the catale of a gros topos, as a sheaf of sheaves, is generated from the catale of underlying spaces as a large idempotent geometric action of the underlying space of spaces: the sheafification.

%Geometric morphisms? Idempotents? Lex comonads?

%\section{From equivalences to bijunctions} \label{Sec:bijun}
%%\input{Z60-bijun}
%
%\section{2-Catales: Functors and naturalities without objects}\label{Sec:adj}
%%\input{Z60-twocatl}

%\section{Points of sobriety and objectivity}\label{Sec:euclid}
%\input{Z80-euclid}
%
%
%

\section{Rocks and starlings}\label{Sec:starlings}
% !TEX root = Z0-catale.tex

A catale is a category where the black box enclosing each object is made transparent. The morphisms refract through the boundaries. Some of the boxes contain copies of other boxes, their retracts, and some morphisms refract through several boundaries. All boxes can be recognized by the refractions and reconstructed as from the maximal idempotents fixing them. A catale is thus a cloud of arrows, akin to the vector field of a dynamical system, but subdivided by walls of light and vortices of wind. A catale displays the dynamical system that becomes visible when the objects of a category relinquish their identities, blend with their neighborhoods, and wave together.
\bear
\frac{particles}{waves}\ \  \approx \ \frac{categories}{catales}\ \  \approx \ \ \frac{rocks}{starlings}
\eear
%
%\begin{figure}[!ht]

\begin{minipage}{.45\linewidth}
\begin{center}
\includegraphics[height=5cm]{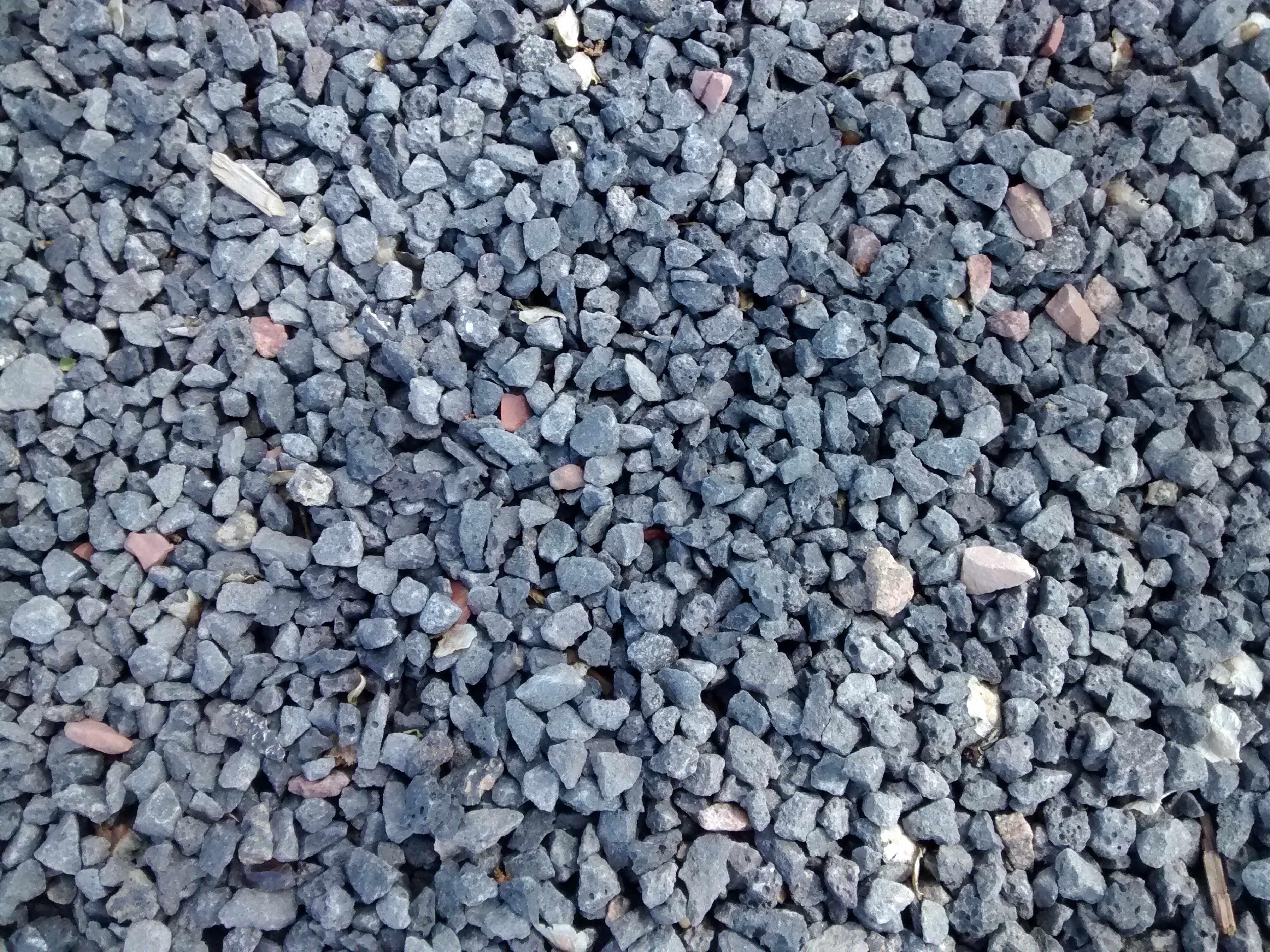}
\end{center}
\end{minipage}
\hspace{.05\linewidth}
\begin{minipage}{.4\linewidth}
\begin{center}
\includegraphics[height=5cm]{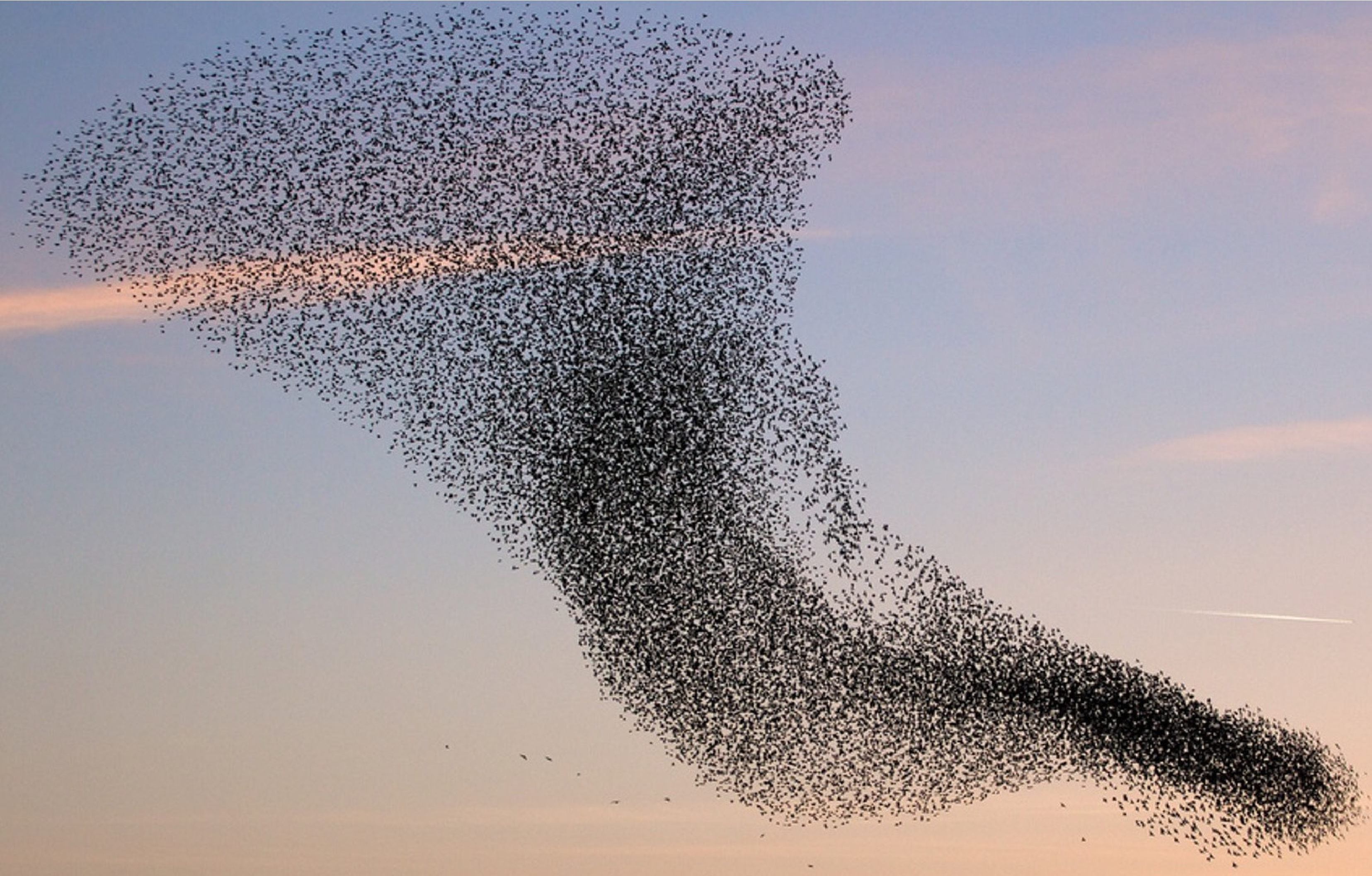}
\end{center}
\end{minipage}
%\caption{Objects have weight. Arrows hide weight.}
%\label{Fig:set-type}
%\end{figure}

\bibliographystyle{plain}
\bibliography{CT,HT,induction,language,PavlovicD,philosophy,ScottDS}

%\appendix
%\input{Z9-Appendix}

\end{document}